\RequirePackage{ifpdf}
\ifpdf 
\documentclass[pdftex]{sigma}
\else
\documentclass{sigma}
\fi
\numberwithin{equation}{section}

\begin{document}

\allowdisplaybreaks

\renewcommand{\PaperNumber}{095}

\FirstPageHeading

\renewcommand{\thefootnote}{$\star$}

\ShortArticleName{Stanilov--Tsankov--Videv Theory}

\ArticleName{Stanilov--Tsankov--Videv Theory\footnote{This paper is a
contribution to the Proceedings of the 2007 Midwest
Geometry Conference in honor of Thomas~P.\ Branson. The full collection is available at
\href{http://www.emis.de/journals/SIGMA/MGC2007.html}{http://www.emis.de/journals/SIGMA/MGC2007.html}}}

\Author{Miguel BROZOS-V\'AZQUEZ~$^{\dag^1}$, Bernd FIEDLER~$^{\dag^2}$,
        Eduardo GARC\'{I}A-R\'IO~$^{\dag^1}$, \\[1mm]
        Peter GILKEY~$^{\dag^3}$, Stana NIK\v CEVI\'C~$^{\dag^4}$,
        Grozio STANILOV~$^{\dag^5}$,\\[1mm]
         Yulian TSANKOV~$^{\dag^5}$, Ram\'on V\'AZQUEZ-LORENZO~$^{\dag^1}$ and Veselin VIDEV~$^{\dag^6}$}

\AuthorNameForHeading{M.~Brozos-V\'azquez et al.}

\Address{$^{\dag^1}$~Department of Geometry and Topology, Faculty of Mathematics,\\
$\phantom{^{\dag^1}}$~University of
      Santiago de Compostela, Santiago de Compostela 15782, Spain}
\EmailDD{\href{mailto:mbrozos@usc.es}{mbrozos@usc.es}, \href{mailto:xtedugr@usc.es}{xtedugr@usc.es}, \href{mailto:ravazlor@usc.es}{ravazlor@usc.es}}

\Address{$^{\dag^2}$~Eichelbaumstr. 13, D-04249 Leipzig, Germany}
\EmailDD{\href{mailto:bfiedler@fiemath.de}{bfiedler@fiemath.de}}

\Address{$^{\dag^3}$~Mathematics Department, University of Oregon, Eugene Oregon 97403-1222, USA}
\EmailDD{\href{mailto:gilkey@uoregon.edu}{gilkey@uoregon.edu}}

\Address{$^{\dag^4}$~Mathematical Institute, SANU, Knez Mihailova 35, p.p. 367,
      11001 Belgrade, Serbia}
\EmailDD{\href{mailto:stanan@mi.sanu.ac.yu}{stanan@mi.sanu.ac.yu}}

\Address{$^{\dag^5}$~Sofia University ``St. Kl. Ohridski'', Sofia, Bulgaria}
\EmailDD{\href{mailto:stanilov@fmi.uni-sofia.bg}{stanilov@fmi.uni-sofia.bg},  \href{mailto:ucankov@fmi.uni-sofia.bg}{ucankov@fmi.uni-sofia.bg}}

\Address{$^{\dag^6}$~Mathematics Department, Thracian University, University Campus,\\
$\phantom{^{\dag^6}}$~6000 Stara Zagora, Bulgaria}
\EmailDD{\href{mailto:videv@uni-sz.bg}{videv@uni-sz.bg}}

\ArticleDates{Received August 07, 2007, in f\/inal form September 22,
2007; Published online September 28, 2007}

\Abstract{We survey some recent results concerning Stanilov--Tsankov--Videv theory, conformal Osserman geometry,
and Walker geometry which relate algebraic properties of the curvature operator to the underlying geometry of
the manifold.}


\Classification{53B20}

\def\MM{\mathfrak{M}}

\begin{flushright}\begin{minipage}{14cm}
\it This article is dedicated to the memory of N.~Bla{\v z}i{\'c} (who passed away 10 October 2005)
and to the memory of T.~Branson (who passed away 11 March 2006). They were coauthors, friends, and talented mathematicians.
\end{minipage}
\end{flushright}

\section{Introduction}
In this article we shall survey
just a few of the many recent developments in Dif\/ferential Geometry which relate algebraic properties of various operators naturally
associated with the curvature of a pseudo-Riemannian manifold to the underlying geometric properties of the manifolds involved.

We introduce the following notational conventions.
Let $\mathcal{M}=(M,g)$ be a pseudo-Riemannian manifold of signature $(p,q)$ and dimension $m=p+q$. We say that $\mathcal{M}$ is
Riemannian if
$p=0$, i.e.\ if $g$ is positive def\/inite. We say that $\mathcal{M}$ is Lorentzian if $p=1$. Let
\[
S^\pm_P(\mathcal{M})=\{\xi\in T_PM:g(\xi,\xi)=\pm1\}
\]
be the {\it pseudo-spheres} of unit spacelike $(+)$ and unit timelike $(-)$ vectors.
Let $\nabla$ be the Levi-Civita connection and let
\[
\mathcal{R}(x,y):=\nabla_x\nabla_y-\nabla_y\nabla_x-\nabla_{[x,y]}
\]
be the associated {\it skew-symmetric curvature operator}. If $\{e_i\}$ is a local frame for the tangent bundle, we let
$g_{ij}:=g(e_i,e_j)$ and let $g^{ij}$ be the inverse matrix. The {\it Jacobi operator} and the {\it Ricci operator} are the self-adjoint
endomorphisms def\/ined, respectively, by:
\begin{gather}\label{eqn-1.a}
\mathcal{J}(x):y\rightarrow\mathcal{R}(y,x)x\qquad\text{and}\qquad\rho:x\rightarrow\sum_{ij}g^{ij}\mathcal{R}(x,e_i)e_j.
\end{gather}
One also def\/ines the {\it curvature tensor} $R\in\otimes^4T^*M$, the {\it scalar curvature} $\tau$, the
{\it Weyl conformal curvature operator} $\mathcal{W}$, and the {\it conformal Jacobi operator} $\mathcal{J}_W$,
respectively, by:
\begin{gather}
R(x,y,z,w)=g(\mathcal{R}(x,y)z,w),\nonumber\\
\tau:=\operatorname{Tr}(\rho)=\textstyle\sum_{ijkl}g^{il}g^{jk}R(e_i,e_j,e_k,e_l),\label{eqn-1.b}\\
\mathcal{W}(x,y): \ \ z\rightarrow\mathcal{R}(x,y)z-\{(m-1)(m-2)\}^{-1}\tau\{g(y,z)x-g(x,z)y\}\nonumber\\
\phantom{\mathcal{W}(x,y): \ \ }{} +(m-2)^{-1}\left\{g(\rho y,z)x-g(\rho x,z)y+g(y,z)\rho x-g(x,z)\rho y\right\},\nonumber\\
\mathcal{J}_W(x): \ \ y\rightarrow\mathcal{W}(y,x)x.\nonumber
\end{gather}

Motivated by the seminal paper of Osserman \cite{Os80}, one studies the spectral properties of the Jacobi operator $\mathcal{J}$ and of
the conformal Jacobi operator $\mathcal{J}_W$ and makes the following def\/initions:

\begin{definition}\label{defn-1}
\rm Let $\mathcal{M}$ be a pseudo-Riemannian manifold.
\begin{enumerate}\itemsep=0pt
\item $\mathcal{M}$ is {\it pointwise Osserman} if $\mathcal{J}$ has constant eigenvalues on
$S^+_P(\mathcal{M})$ and on $S^-_P(\mathcal{M})$ for every $P\in M$.
\item $\mathcal{M}$ is {\it pointwise conformally Osserman} if $\mathcal{J}_W$ has constant eigenvalues on
$S^+_P(\mathcal{M})$ and on~$S^-_P(\mathcal{M})$ for every $P\in M$.
\end{enumerate}\end{definition}

We refer to \cite{GaKuVa02} for a more complete discussion of Osserman geometry as that lies beyond the scope of our
present endeavors.

Similarly, motivated by the seminal papers of Stanilov and Videv \cite{StVi04}, of Tsankov \cite{Ts05}, and of Videv \cite{Vixx} one
studies the commutativity properties of these operators:

\begin{definition}\label{defn-2}
\rm Let $\mathcal{M}$ be a pseudo-Riemannian manifold.
\begin{enumerate}\itemsep=0pt
\item $\mathcal{M}$ is {\it Jacobi--Tsankov} if $\mathcal{J}(\xi_1)\mathcal{J}(\xi_2)=\mathcal{J}(\xi_2)\mathcal{J}(\xi_1)$
for all $\xi_i$.
\item $\mathcal{M}$ is {\it mixed-Tsankov} if $\mathcal{R}(\xi_1,\xi_2)\mathcal{J}(\xi_3)=\mathcal{J}(\xi_3)\mathcal{R}(\xi_1,\xi_2)$
for all $\xi_i$.
\item $\mathcal{M}$ is {\it skew-Tsankov} if $\mathcal{R}(\xi_1,\xi_2)\mathcal{R}(\xi_3,\xi_4)
    =\mathcal{R}(\xi_3,\xi_4)\mathcal{R}(\xi_1,\xi_2)$ for all $\xi_i$.
\item $\mathcal{M}$ is {\it Jacobi--Videv} if $\mathcal{J}(\xi)\rho=\rho \mathcal{J}(\xi)$
for all $\xi$.
\item $\mathcal{M}$ is {\it skew-Videv} if $\mathcal{R}(\xi_1,\xi_2)\rho=\rho\mathcal{R}(\xi_1,\xi_2)$ for all $\xi_i$. This has also
been called {\it {\it Ricci semi-symmetric}} by some authors.
\end{enumerate}\end{definition}

In this brief note, we survey some recent results concerning these concepts; we refer to~\cite{GaKuVa02,Gi01,Gi07a} for a discussion of
some previous results in this area.

Our f\/irst task is to pass to the algebraic setting.
\begin{definition}\label{defn-3}
Let $\langle\cdot,\cdot\rangle$ be a non-degenerate bilinear form of signature $(p,q)$ on a f\/inite dimensional real vector space $V$. Let
$R\in\otimes^4V^*$ be a $4$-tensor. We say that
$\MM=(V,\langle\cdot,\cdot\rangle,R)$ is a~{\it model} and that $R$ is an {\it algebraic curvature tensor} if $R$ satisf\/ies the usual
curvature identities for all
$x$,~$y$,~$z$, and~$w$:
\begin{gather*}
R(x,y,z,w)=-R(y,x,z,w)=R(z,w,x,y),\\
R(x,y,z,w)+R(y,z,x,w)+R(z,x,y,w)=0.
\end{gather*}
The associated {\it algebraic curvature operator} $\mathcal{R}$ is then def\/ined by using the inner product to raise indices; this
skew-symmetric operator is characterized by the identity:
\[
\langle\mathcal{R}(x,y)z,w\rangle=R(x,y,z,w).
\]
The Jacobi operator, the Ricci operator, the Weyl conformal curvature operator, and the conformal Jacobi operator are then def\/ined as in
equations (\ref{eqn-1.a}) and (\ref{eqn-1.b}). The concepts of Def\/initions~\ref{defn-1} and~\ref{defn-2} extend naturally to this
setting.\end{definition}

 If $P$ is a point of a pseudo-Riemannian manifold $\mathcal{M}$, then the
associated model is def\/ined by
\[
\MM(\mathcal{M},P):=(T_PM,g_P,R_P).
\]
 We note that every model $\MM$
is geometrically realizable; this means that given $\MM$, there is $(\mathcal{M},P)$ such that
$\MM(\mathcal{M},P)$ is isomorphic to $\MM$ -- see, for example, the discussion in \cite{Gi01}.

One has the following examples of algebraic curvature tensors.

\begin{example}\label{exm-1}
\ \begin{enumerate}\itemsep=0pt \item
If $\psi$ is self-adjoint with respect to $\langle\cdot,\cdot\rangle$, one def\/ines an algebraic curvature tensor
\[
R_\psi(x,y,z,w)=\langle\psi x,w\rangle\langle\psi y,z\rangle-\langle\psi x,z\rangle\langle\psi y,w\rangle.
\]
Taking $\psi=\operatorname{id}$ and rescaling yields the algebraic curvature tensor of constant sectional curvature
$c$:
\[
R_c(x,y,z,w)=c\{\langle x,w\rangle\langle y,z\rangle-\langle x,z\rangle\langle y,w\rangle\}.
\]
One says that a model $\MM$ or a pseudo-Riemannian manifold $\mathcal{M}$ has {\it constant sectional curvature $c$} if $R=R_c$ for
some constant $c$.
\item If $\phi$ is skew-adjoint with respect to $\langle\cdot,\cdot\rangle$, one def\/ines an algebraic curvature tensor
\[
R_\phi(x,y,z,w)=\langle\phi y,z\rangle\langle\phi x,w\rangle-\langle\phi x,z\rangle\langle\phi y,w\rangle
-2\langle\phi x,y\rangle\langle\phi z,w\rangle.
\]
\end{enumerate}\end{example}

\begin{remark}\label{rmk-1}
 The space of algebraic curvature tensors is spanned as a linear space by the tensors given in
Example \ref{exm-1} (1) or in Example~\ref{exm-1} (2) \cite{Fi03}; we also refer to \cite{DiGa04}.
\end{remark}

Our f\/irst result is the equivalence of conditions (1) and (2) and of (4) and (5) in Def\/inition~\ref{defn-2}; if $\MM$ is a model or if
$\mathcal{M}$ is a pseudo-Riemannian manifold, then Jacobi--Tsankov and mixed-Tsankov are equivalent conditions. Similarly
Jacobi--Videv and skew-Videv are equivalent conditions. This follows from the following result \cite{GiNi07}:
\begin{theorem}\label{thm-1} Let $\MM$ be a model and let $T$ be a self-adjoint linear transformation of $V$. Then
the following assertions are equivalent:
\begin{enumerate}\itemsep=0pt
\item[\rm 1.] $\mathcal{R}(x,y)T=T\mathcal{R}(x,y)$ for all $x$, $y\in V$.
\item[\rm 2.] $\mathcal{J}(x)T=T\mathcal{J}(x)$ for all $x\in V$.
\item[\rm 3.] $R(Tx,y,z,w)=R(x,Ty,z,w)=R(x,y,Tz,w)=R(x,y,z,Tw)$ for all $x$, $y$, $z$, $w$ in $V$.
\end{enumerate}\end{theorem}

Here is a brief outline of the remainder of this article. In Section \ref{sect-2}, we study Jacobi--Tsankov models and
manifolds. In Section \ref{sect-3}, we study skew-Tsankov models and manifolds. In Section~\ref{sect-4}, we study
Jacobi--Videv models and manifolds. In Section \ref{sect-5}, we recall some general results concerning conformal Osserman
geometry. In Section~\ref{sect-6}, we study these concepts in the context of Walker manifolds of signature (2,2).

\section[Jacobi-Tsankov models and manifolds]{Jacobi--Tsankov models and manifolds}\label{sect-2}

We f\/irst turn to the Riemannian setting in the following result \cite{BrGi07}:

\begin{theorem}\label{thm-2}
If $\MM$ is a Jacobi--Tsankov Riemannian model, then $R=0$.
\end{theorem}

\begin{proof} We can sketch the proof as follows. Since $\{\mathcal{J}(x)\}_{x\in V}$ form a family of
commuting self-adjoint operators, we can simultaneously diagonalize these operators to
decompose $V=\oplus_\lambda V_\lambda$
so $\mathcal{J}(x)=\lambda(x)\operatorname{id}$ on $V_\lambda$. If $x\in V$, decompose $x=\oplus x_\lambda$ for $x_\lambda\in
V_\lambda$. Let
\[
\mathcal{O}=\{x\in V:x_\lambda\ne0\quad\text{for all}\quad\lambda\};
\]
this is an open dense
subset of $V$. If $x\in\mathcal{O}$, since $\mathcal{J}(x)x=0$, $\lambda(x)=0$ for all
$\lambda$. Since $\mathcal{O}$ is dense and $\lambda(\cdot)$ is continuous,
$\lambda(x)=0$ for all $x$ so $\mathcal{J}(x)=0$ for all $x$; the usual curvature symmetries now imply
the full curvature tensor $R$ vanishes.\end{proof}

\begin{definition}\label{defn-4}
One says that a model
$\MM$ or a pseudo-Riemannian manifold
$\mathcal{M}$ is {\it orthogonally Jacobi--Tsankov} if
$\mathcal{J}(x)\mathcal{J}(y)=\mathcal{J}(y)\mathcal{J}(x)$ for all vectors $x$ and $y$ with $x\perp y$.
\end{definition}

One has the following
classif\/ication result
\cite{BrGi07}; we also refer to a related result
\cite{Ts05} if
$\mathcal{M}$ is a~hypersurface in $\mathbb{R}^{m+1}$.

\begin{theorem}\label{thm-3}
\ \begin{enumerate}\itemsep=0pt
\item[\rm 1.] Let $\MM=(V,\langle\cdot,\cdot\rangle,R)$ be a Riemannian model. Then $\MM$
is orthogonally Jacobi--Tsankov if and only if one of the following conditions holds:
\begin{enumerate}\itemsep=0pt
\item[\rm (a)] $R=cR_{\operatorname{id}}$ has constant sectional curvature $c$ for some $c\in\mathbb{R}$.
\item[\rm (b)] $\dim(V)$ is even and $R=cR_\Theta$ is defined by Example {\rm \ref{exm-1} (2)} where $\Theta$ is a Hermitian
almost complex structure on
$(V,\langle\cdot,\cdot\rangle)$ and where $c\in\mathbb{R}$.
\end{enumerate}\itemsep=0pt
\item[\rm 2.] Let $\mathcal{M}$ be a Riemannian manifold of dimension $m$.\begin{enumerate}
\item[\rm (a)] If $m>2$, then $\mathcal{M}$ is orthogonally
Jacobi--Tsankov if and only if
$\mathcal{M}$ has constant sectional curvature $c$.
\item[\rm (b)]
If $m=2$, then $\mathcal{M}$ is always orthogonally Jacobi--Tsankov.
\end{enumerate}\end{enumerate}
\end{theorem}

\begin{definition}\label{defn-5} We say that a model $\MM$ or a pseudo-Riemannian manifold $\mathcal{M}$ is {\it conformally
Jacobi--Tsankov} if
$\mathcal{J}_W(x)\mathcal{J}_W(y)=\mathcal{J}_W(y)\mathcal{J}_W(x)$ for all $x$ and $y$.
We say that $\MM$ or $\mathcal{M}$ is {\it orthogonally conformally Jacobi--Tsankov} if
$\mathcal{J}_W(x)\mathcal{J}_W(y)=\mathcal{J}_W(y)\mathcal{J}_W(x)$ for all vectors~$x$ and~$y$ with $x\perp  y$.
\end{definition}

\begin{remark}\label{rmk-2} These are  conformal notions -- if $\mathcal{M}$
is conformally equivalent to $\mathcal{M}_1$, then $\mathcal{M}$ is conformally Jacobi--Tsankov (resp.\
orthogonally conformally Jacobi--Tsankov) if and only if
$\mathcal{M}_1$ is conformally Jacobi Tsankov (resp.\ orthogonally conformally Jacobi--Tsankov). We refer to~\cite{BlGiNiSi05} for further details.\end{remark}

We have the following useful result:

\begin{theorem}\label{thm-4}
A Riemannian model $\MM$ is orthogonally conformally
Jacobi--Tsankov if and only if $\mathcal{W}=0$.
\end{theorem}

\begin{proof} Let $\mathcal{W}$ be the associated Weyl conformal curvature operator. Then $\mathcal{W}$ is an
algebraic curvature tensor which is orthogonally-Jacobi Tsankov. Thus Theorem~\ref{thm-3}
yields either that $\mathcal{W}=c\mathcal{R}_{\operatorname{id}}$ or that $\mathcal{W}=c\mathcal{R}_\Theta$. Since the scalar
curvature def\/ined by the tensors $\mathcal{R}_{\operatorname{id}}$ and $\mathcal{R}_\Theta$ is non-zero, we may
conclude
$c=0$.
\end{proof}

There are non-trivial examples of Jacobi--Tsankov manifolds and models in the higher signature setting.
\begin{definition}\label{defn-6} We say that a model
$\MM=(V,\langle\cdot,\cdot\rangle,R)$ is {\it indecomposable} if there is no non-trivial orthogonal decomposition
$V=V_1\oplus V_2$ which induces a decomposition $R=R_1\oplus R_2$.\end{definition}

 We refer to \cite{BrGi06} for the proof of the following result:

\begin{theorem}\label{thm-5}
Let $\MM$ be a model.
\begin{enumerate}
\item[\rm 1.] If $\MM$ is Jacobi--Tsankov, then $\mathcal{J}(x)^2=0$ for all $x$ in $V$.
\item[\rm 2.] If $\MM$ is Jacobi--Tsankov and Lorentzian, then $R=0$.
\item[\rm 3.] Let $\MM$ be indecomposable with $\dim(\MM)<14$. The following conditions are equivalent:
\begin{enumerate}\itemsep=0pt
\item[\rm (a)] $V=U\oplus\bar U$ and $R=R_U\oplus0$ where $U$ and $\bar U$ are totally isotropic subspaces.
\item[\rm (b)] $\MM$ is Jacobi--Tsankov.
\end{enumerate}
Either {\rm (3a)} or {\rm (3b)} implies that $\mathcal{R}(x,y)z\in\bar U$ and that $\mathcal{R}(x,y)\mathcal{R}(u,v)z=0$ for all
$x$, $y$, $z$, $u$, $v\in V$, that
$\mathcal{J}(x)\mathcal{J}(y)=0$ for all $x$, $y\in V$, and that $\MM$ is skew-Tsankov.
\end{enumerate}\end{theorem}

The condition $\mathcal{J}(x)^2=0$ does not imply that $\MM$ is Jacobi--Tsankov \cite{BrGi06}:
\begin{example}\label{exm-2}
Let $\langle\cdot,\cdot\rangle$ be an inner product of signature $(4,4)$ on $\mathbb{R}^8$.
Choose skew-symmetric endomorphisms $\{e_1,e_2,e_3,e_4\}$ so that
\[
e_1^2=e_2^2=\operatorname{id},\qquad e_3^2=e_4^2=-\operatorname{id},
\qquad\text{and}\qquad e_ie_j+e_je_i=0\quad\text{for}\quad i\ne j.
\]
Note that this gives a suitable
Clif\/ford module structure to
$\mathbb{R}^8$. Set
\[
\phi_1=e_1+e_3\qquad\text{and}\qquad\phi_2=e_2+e_4.
\]
 Adopt the notation of Example \ref{exm-1} (2) to def\/ine
$R_{\phi_i}$. Then
\[
\MM:=(\mathbb{R}^8,\langle\cdot,\cdot\rangle,R_{\phi_1}+R_{\phi_2})
\] is not Jacobi--Tsankov but
satisf\/ies $\mathcal{J}(x)^2=0$ for all
$x$.\end{example}

We have the following example \cite{BrGi06} that shows that the structure of Theorem~\ref{thm-5} (3a) is geometrically realizable:

\begin{example}\label{exm-3}
Let $(x_1,\dots,x_p,y_1,\dots,y_p)$ be the usual coordinates
on $\mathbb{R}^{2p}$. Let $\mathcal{M}=(\mathbb{R}^{2p},g)$ where
$g(\partial_{x_i},\partial_{y_j})=\delta_{ij}$ and let
$g(\partial_{x_i},\partial_{x_j})=g_{ij}(x)$.
Then there exists a decomposition $T(\mathbb{R}^{2p})=U\oplus\bar U$ where $U$ and
$\bar U$ are totally isotropic so that $\mathcal{R}(x,y)z\in\bar U$ and that $\mathcal{R}(x,y)\mathcal{R}(u,v)z=0$ for all
$x$, $y$, $z$, $u$, $v\in V$. Furthermore, for
generic
$g$, the model $\MM(\mathcal{M},P)$ is indecomposable for all $P\in\mathbb{R}^{2p}$.
\end{example}

The restriction in Theorem \ref{thm-5} that $\dim(V)<14$ is essential. We have the following \cite{BrGi06}:

\begin{example}\label{exm-4}
Let $\{\alpha_i,\alpha_i^*,\beta_{i,1},\beta_{i,2},\beta_{4,1},\beta_{4,2}\}_{1\le i\le 3}$ be a basis
for $\mathbb{R}^{14}$. Def\/ine $\mathcal{M}_{6,8}$ by:
\begin{gather*}
\langle\alpha_i,\alpha_i^*\rangle=\langle\beta_{i,1},\beta_{i,2}\rangle,\quad 1\le i\le 3;
\qquad \langle\beta_{4,1},\beta_{4,1}\rangle=\langle\beta_{4,2},\beta_{4,2}\rangle=-\tfrac12;\qquad
\langle\beta_{4,1},\beta_{4,2}\rangle=\tfrac14;\\
R_{\alpha_2,\alpha_1,\alpha_1,\beta_{2,1}}
=R_{\alpha_3,\alpha_1,\alpha_1,\beta_{3,1}}
=R_{\alpha_3,\alpha_2,\alpha_2,\beta_{3,2}}=1,\\
R_{\alpha_1,\alpha_2,\alpha_2,\beta_{1,2}}
=R_{\alpha_1,\alpha_3,\alpha_3,\beta_{1,1}}
=R_{\alpha_2,\alpha_3,\alpha_3,\beta_{2,2}}=1,\\
R_{\alpha_1,\alpha_2,\alpha_3,\beta_{4,1}}
=R_{\alpha_1,\alpha_3,\alpha_2,\beta_{4,1}}
=R_{\alpha_2,\alpha_3,\alpha_1,\beta_{4,2}}
=R_{\alpha_2,\alpha_1,\alpha_3,\beta_{4,2}}=-\tfrac{1}{2}.
\end{gather*}
Then $\MM_{6,8}$ has signature $(6,8)$, $\MM_{6,8}$ is Jacobi--Tsankov, $\MM_{6,8}$ is not
skew-Tsankov, and there exist $x$ and $y$ so that
$\mathcal{J}(x)\mathcal{J}(y)\ne0$.\end{example}

Furthermore, this example is geometrically realizable \cite{BrGiNi07}:

\begin{example}\label{exm-5}
Take coordinates
$\{x_i,x_i^*,y_{i,1},y_{i,2},y_{4,1},y_{4,2}\}_{i=1,2,3}$
for $\mathbb{R}^{14}$. Let $a_{i,j}\in\mathbb{R}$ and let
$\mathcal{M}_{6,8}:=(\mathbb{R}^{14},g)$ where:
\begin{gather*}
{g}(\partial_{x_i},\partial_{x_i^*})=g(\partial_{y_{i,1}},\partial_{y_{i,2}})=1,\qquad
{g}(\partial_{y_{4,1}},\partial_{y_{4,1}})=g(\partial_{y_{4,2}},\partial_{y_{4,2}})=-\tfrac12,\\
{g}(\partial_{y_{4,1}},\partial_{y_{4,2}})=\tfrac14,\qquad
{g}(\partial_{x_1},\partial_{x_1})=-2 a_{2,1} x_2 y_{2,1}-2 a_{3,1}x_3 y_{3,1},\\
{g}(\partial_{x_2},\partial_{x_2})=-2 a_{3,2}x_3 y_{3,2}-2 a_{1,2}x_1 y_{1,2},\qquad
{g}(\partial_{x_3},\partial_{x_3})=-2 a_{1,1}x_1 y_{1,1}-2 a_{2,2}x_2 y_{2,2},\\
{g}(\partial_{x_1},\partial_{x_2})=2 (1-a_{2,1}) x_1
y_{2,1}+2(1-a_{1,2}) x_2 y_{1,2},\\
{g}(\partial_{x_2},\partial_{x_3})=x_1y_{4,1}+2 (1-a_{3,2})x_2 y_{3,2}+2(1-a_{2,2})x_3 y_{2,2},\\
{g}(\partial_{x_1},\partial_{x_3})=x_2y_{4,2}+2(1-a_{3,1})x_1
y_{3,1}+2(1-a_{1,1})x_3y_{1,1}.
\end{gather*}
Then $\mathcal{M}$ has the model $\MM_{6,8}$ and $\mathcal{M}$ is locally symmetric if and only if
\begin{gather*}
a_{1,1}+a_{2,2}+a_{3,1}a_{3,2}=2,\qquad
3 a_{2,1}+3a_{3,1}+3 a_{1,2} a_{1,1}=4,\\
3 a_{1,2}+3a_{3,2}+3 a_{2,1} a_{2,2}=4.
\end{gather*}
\end{example}

We note that the relations of Example \ref{exm-5} have non-trivial solutions. One may take, for example,
$a_{1,1}=a_{2,2}=1$, $a_{1,2}=a_{2,1}=\frac23$, and  $a_{3,1}=a_{3,2}=0$.

\section{Skew-Tsankov models and manifolds}\label{sect-3}

Riemannian skew-Tsankov models are completely classif\/ied
\cite{BrGi07b}:

\begin{theorem}\label{thm-6}
Let $\MM$ be a Riemannian skew-Tsankov
model. Then there is an orthogonal direct sum decomposition $V=V_1\oplus\cdots\oplus V_k\oplus U$ where
$\dim(V_k)=2$ and where
$R=R_1\oplus\cdots\oplus R_k\oplus 0$.
\end{theorem}

\begin{proof}  One has that $\{\mathcal{R}(\xi,\eta)\}_{\xi,\eta\in V}$ is a collection of commuting
skew-adjoint endomorphisms. As the inner product is def\/inite, there exists an orthogonal decomposition
of $V$ so that each endomorphism $\mathcal{R}(\xi,\eta)$ decomposes as a direct sum of $2\times 2$ blocks
\[
\left(\begin{array}{rr}
    0&a(\xi,\eta)\\-a(\xi,\eta)&0
   \end{array}\right).
\]
The desired result then follows from the curvature symmetries.
\end{proof}

The situation in the geometric context is less clear. We refer to \cite{BrGi07b} for the following $3$-di\-mensional and $4$-dimensional
examples which generalize previous examples found in \cite{Ts05}. We say that $\mathcal{M}$ is an irreducible Riemannian manifold if
there is no local product decomposition.

\begin{example}\label{exm-6}
\ \begin{enumerate}\itemsep=0pt \item Let $M=(0,\infty)\times N$ where $N$ is a Riemann surface
with scalar curvature $\tau_N\not\equiv1$. Give $M$ the warped product metric $ds^2=dt^2+t^2ds^2_N$. Then
$\mathcal{M}:=(M,g_M)$ is an irreducible skew-Tsankov manifold with $\tau_M=t^{-2}(\tau_N-1)$.
\item Let  $(x_1,x_2,x_3,x_4)$ be the usual coordinates on $\mathbb{R}^4$. Let $\mathcal{M}_\beta=(\mathbb{R}^4,g)$
where
$ds^2=x_3^2dx_1^2+(x_3+\beta x_4)^2dx_2^2+dx_3^2+dx_4^2$. Then $\mathcal{M}_\beta$ is an
irreducible skew-Tsankov manifold with $\tau_\beta=-2x_3^{-1}(x_3+\beta
x_4)^{-1}$. $M_\beta$ is not isometric to $M_{\bar\beta}$ for
$0<\beta<\bar\beta$.\end{enumerate}
\end{example}

In the higher signature setting, we note that Example \ref{exm-3} provides examples of neutral signature pseudo-Riemannian manifolds
with $\mathcal{R}(x,y)\mathcal{R}(z,w)=0$ for all $x$, $y$, $z$, $w$. There are, however, less trivial examples.
\begin{definition}\label{defn-7} We say $\mathcal{M}$ is $3$-skew nilpotent if
\begin{enumerate}\itemsep=0pt
\item There exist $\xi_i$ with $\mathcal{R}(\xi_1,\xi_2)\mathcal{R}(\xi_3,\xi_4)\ne0$ and
\item For all $\xi_i$, one has $\mathcal{R}(\xi_1,\xi_2)\mathcal{R}(\xi_3,\xi_4)\mathcal{R}(\xi_5,\xi_6)=0$.
\end{enumerate}\end{definition}
We refer to \cite{FiGi03} for the
proof of:

\begin{example}\label{exm-7}
Let $(x,u_1,\dots,u_{m-2},y)$ be coordinates on $\mathbb{R}^m$. Let $f=f(\vec u)$ be smooth.
Let
$\Xi$ be a~non-degenerate bilinear form on
$\mathbb{R}^{m-2}$. Consider $\mathcal{M}:=(\mathbb{R}^m,g)$ where the non-zero components of $g$ are given by:
\[
g(\partial_x,\partial_x)=-2f(\vec u),\qquad
g(\partial_x,\partial_y)=1,\qquad
g(\partial_{u_a},\partial_{u_b})=\Xi_{ab}.
\]
Then $\mathcal{M}$ is skew-Tsankov and
$3$-skew nilpotent; it need not be Jacobi--Tsankov.
\end{example}

\section[Jacobi-Videv models and manifolds]{Jacobi--Videv models and manifolds}\label{sect-4}

One says $\mathcal{M}$ is Einstein if $\rho$ is a scalar multiple of the identity. More generally:
\begin{definition}\label{defn-8} One
says
$\mathcal{M}$ is {\it pseudo-Einstein} if $\rho$ either has a single real eigenvalue $\lambda$ or has exactly two
eigenvalues which are complex conjugates $\mu$ and $\bar\mu$.\end{definition}
It is immediate that pseudo-Einstein implies Einstein in the Riemannian setting as $\rho$ is diagonalizable if the metric
is positive def\/inite.

 We refer to \cite{GiPuVi07} for the proof of the following
result; see also \cite{IvViZh07} for related work in the $4$-dimensional context.

\begin{theorem}\label{thm-7}
Let $\MM$ be an indecomposable model which is Jacobi--Videv. Then
$\MM$ is pseudo-Einstein.
\end{theorem}

\begin{proof} Let $m:=\dim(V)$. Let $\lambda\in\mathbb{C}$ have non-negative real part. Set
\[
V_\lambda:=\{v\in V:(T-\lambda)^m(T-\bar\lambda)^mv=0\}.
\]
We then have the Jordan decomposition $V=\oplus_\lambda V_\lambda$ as an orthogonal direct sum of generalized
eigenspaces of $\rho$. Since
$\mathcal{J}(x)$ preserves this decomposition, it follows that $\mathcal{J}=\oplus_\lambda\mathcal{J}_\lambda$.
The curvature symmetries then imply that $\mathcal{R}=\oplus_\lambda\mathcal{R}_\lambda$. Since $\MM$ is assumed indecomposable,
there is only one $V_\lambda\ne\{0\}$ and thus $\MM$ is pseudo-Einstein.
\end{proof}

This shows, in the Riemannian setting, that an indecomposable model is Jacobi--Videv if and only if it is Einstein.
The condition that $\MM$ is pseudo-Einstein does not, however, imply that $\MM$ is Jacobi--Videv in the higher signature setting
as the following \cite{GiNi07} shows:

\begin{example}\label{exm-8}
Let $\{x_1,x_2,x_3,x_4\}$ be coordinates on $\mathbb{R}^4$.
 Let $\mathcal{M}=(\mathbb{R}^4,g)$ where
\[
g(\partial_{x_1},\partial_{x_4})=g(\partial_{x_2},\partial_{x_2})=g(\partial_{x_3},\partial_{x_3})=1\qquad\text{and}\qquad
g(\partial_{x_1},\partial_{x_3})=e^{x_2}.
\]
Then $\mathcal{M}$ is a
homogeneous Lorentz manifold and $\mathcal{M}$ is pseudo-Einstein with
$\operatorname{Rank}(\rho)=2$, $\operatorname{Rank}(\rho^2)=1$, and
$\operatorname{Rank}(\rho^3)=0$. Thus $\mathcal{M}$ is pseudo-Einstein.
However $\mathcal{M}$ is not Jacobi--Videv.\end{example}

We also have \cite{GiNi07}
\begin{example}\label{exm-9}
Let $\{x_1,x_2,x_3,x_4\}$ be coordinates on $\mathbb{R}^4$.
 Let $\mathcal{M}=(\mathbb{R}^4,g)$ where
\begin{gather*}
g(\partial_{x_1},\partial_{x_3})=g(\partial_{x_2},\partial_{x_4})=1,\qquad
g(\partial_{x_3},\partial_{x_3})=-g(\partial_{x_4},\partial_{x_4})=sx_1x_2,\\
g(\partial_{x_3},\partial_{x_4})=\tfrac s2(x_2^2-x_1^2).
\end{gather*}
 Then $\mathcal{M}$ is locally symmetric of signature $(2,2)$, $\cal{M}$ is Jacobi--Videv,  $\cal{M}$ is skew-Tsankov, and
$\cal{M}$ is conformal Osserman.
$\mathcal{M}$ is neither Jacobi--Tsankov nor Osserman. $\mathcal{M}$ is pseudo-Einstein with $\rho^2=-s^2\operatorname{id}$.
\end{example}

\begin{example}\label{exm-10}
Setting
\begin{gather*}
g(\partial_{x_1},\partial_{x_3})=g(\partial_{x_2},\partial_{x_4})=1,\qquad
  g(\partial_{x_3},\partial_{x_3})=-g(\partial_{x_4},\partial_{x_4})=\tfrac s2(x_2^2-x_1^2),\\
g(\partial_{x_3},\partial_{x_4})=-sx_1x_2
\end{gather*}
yields a local symmetric space of signature $(2,2)$ which is Einstein. This manifold is Jacobi--Videv and skew-Tsankov. It is neither
Jacobi--Tsankov, Osserman, nor conformal Osserman.
\end{example}

We can give a general ansatz which constructs such examples in the algebraic setting; we do not know if these examples are
geometrically realizable in general:

\begin{example}\label{exm-11}
Let $\MM=(V,(\cdot,\cdot),R)$ be a model.  We complexify and let
\[
U:=V\otimes_{\mathbb{R}}\mathbb{C}.
\]
We extend $(\cdot,\cdot)$ and $R$ to be complex multi-linear. Let $\{e_i\}$ be an orthonormal basis for $V$. Let
$\{e_i^+:=e_i,e_i^-:=\sqrt{-1}e_i\}$ be a basis for the underlying real vector space $U:=V\oplus\sqrt{-1}V$. Let
$\Re$ and $\Im$ denote the real and imaginary parts of a complex number, respectively. It is then immediate that
\[
\langle\cdot,\cdot\rangle:=\Re\{(\cdot,\cdot)\}\qquad\text{and}\qquad
  S(\cdot,\cdot,\cdot,\cdot)=\Im\{R(\cdot,\cdot,\cdot,\cdot)\}
  \]
def\/ine a model $\mathfrak{N}:=(U,\langle\cdot,\cdot\rangle,S)$. One has that the non-zero components of
$\langle\cdot,\cdot\rangle$ are
$\langle e_i^+,e_i^+\rangle=1$ and $\langle e_i^-,e_i^-\rangle=-1$. Thus the metric has neutral signature.
Furthermore, the non-zero components of $S$ are given by:
\begin{gather*}
S(e_i^-,e_j^+,e_k^+,e_l^+)=S(e_i^+,e_j^-,e_k^+,e_l^+)=S(e_i^+,e_j^+,e_k^-,e_l^+)
    \\ \qquad{}=S(e_i^+,e_j^+,e_k^+,e_l^-)=R(e_i,e_j,e_k,e_l),\\
S(e_i^+,e_j^-,e_k^-,e_l^-)=S(e_i^-,e_j^+,e_k^-,e_l^-)=S(e_i^-,e_j^-,e_k^+,e_l^-)
    \\ \qquad{}=S(e_i^-,e_j^-,e_k^-,e_l^+)=-R(e_i,e_j,e_k,e_l).
\end{gather*}
\end{example}

We refer to  \cite{GiNi07} for the proof of the following result:

\begin{theorem}\label{thm-8}
Adopt the notation of Example {\rm \ref{exm-11}}. If $\MM$ is a Riemannian Einstein model with $\rho_{\MM}=s\operatorname{id}$, then
$\mathfrak{N}$ is a Jacobi--Videv pseudo-Einstein neutral signature model with $\rho_{\mathfrak{N}}^2=-4s^2\operatorname{id}$.
\end{theorem}

\begin{definition}\label{defn-9}
Let $\MM=(V,\langle\cdot,\cdot\rangle,R)$ be a model. Let $\{v_1,\dots,v_k\}$ be an orthonormal basis for a
non-degenerate
$k$-plane
$\pi\subset V$. Let $\varepsilon_i:=\langle v_i,v_i\rangle$ be $\pm1$. One def\/ines the {\it higher order Jacobi operator}
by setting:
\[
\mathcal{J}(\pi):=\sum_{i=1}^k\varepsilon_i\mathcal{J}(v_i).
\]\end{definition}

The operator $\mathcal{J}(\pi)$ is independent of the particular orthonormal basis chosen; we refer to
\cite{GiStVi98,StVi92,StVi98} for a further discussion of this operator. If $\pi=V$, then $\mathcal{J}(\pi)=\rho$. If
$\pi=\operatorname{Span}(x)$ where
$x$ is a unit spacelike vector, then
$\mathcal{J}(\pi)=\mathcal{J}(x)$. Thus $\mathcal{J}(\pi)$ can be thought of as interpolating between the Jacobi operator and the Ricci
operator.

\begin{definition}\label{defn-10} Let $\MM$ be a model of signature $(p,q)$. We say that $(r,s)$ is {\it admissible} if and only if
\[
0\le r\le p,\qquad 0\le s\le q,\qquad\text{and}\qquad
1\le r+s\le m-1.
\] Equivalently, $(r,s)$ is {\it admissible} if and only if the Grassmannian
of linear subspaces of signature $(r,s)$ has positive dimension.
\end{definition}

 One has the following
useful characterization \cite{GiPuVi07}:
\begin{theorem}\label{thm-9}
\ The following properties are equivalent for $\MM=(V,\langle\cdot,\cdot\rangle,R)$:
\begin{enumerate}\itemsep=0pt
\item[\rm 1.] $\MM$ is Jacobi--Videv, i.e.\ $\mathcal{J}(x)\rho=\rho\mathcal{J}(x)$ for all $x\in V$.
\item[\rm 2.] There exists $(r,s)$ admissible so $\mathcal{J}(\pi)\mathcal{J}(\pi^\perp)=\mathcal{J}(\pi^\perp)\mathcal{J}(\pi)$ for
every non-degenerate subspace $\pi$ of signature $(r,s)$.
\item[\rm 3.] There exists $(r,s)$ admissible so $\mathcal{J}(\pi)\rho=\rho\mathcal{J}(\pi)$ for every
non degenerate subspace $\pi$ of signature $(r,s)$.
\item[\rm 4.] $\mathcal{J}(\pi)\mathcal{J}(\pi^\perp)=\mathcal{J}(\pi^\perp)\mathcal{J}(\pi)$ for every non-degenerate linear subspace
$\pi$.
\item[\rm 5.] $\mathcal{J}(\pi)\rho=\rho\mathcal{J}(\pi)$ for every non-degenerate linear subspace $\pi\subset V$.
\end{enumerate}
\end{theorem}

\section{Conformal Osserman geometry}\label{sect-5}

We refer to \cite{BlGi04a,BlGiNiSi05} for the proof of the following result:

\begin{theorem}\label{thm-10}
Let $\mathcal{M}$ be a conformally Osserman pseudo-Riemannian manifold of dimension~$m$.
\begin{enumerate}\itemsep=0pt \item[\rm 1.] If $\mathcal{M}$ is Riemannian and if $m$ is odd, then $\mathcal{M}$ is locally
conformally flat.
\item[\rm 2.] If $\mathcal{M}$ is Riemannian, if $m\equiv2$ {\rm mod} $4$, if $m\ge10$, and if $\mathcal{W}(P)\ne0$,
then there is an open neighborhood of $P$ in
$M$ which is conformally equivalent to an open subset of either complex projective space with the Fubini--Study
metric or the negative curvature dual.
\item[\rm 3.] If $\mathcal{M}$ is Lorentzian, then $\mathcal{M}$ is locally conformally flat.
\end{enumerate}
\end{theorem}

We also recall the following result \cite{BlGi05,BrGaVa06}:
\begin{theorem}\label{thm-11}
Let $\mathcal{M}$ be a $4$-dimensional model of arbitrary signature.
\begin{enumerate}\itemsep=0pt \item[\rm 1.] $\mathfrak{M}$ is conformally Osserman if and only if
$\mathcal{M}$ is either self-dual or anti-self-dual.
\item[\rm 2.] If $\mathfrak{M}$ is Riemannian, then $\mathfrak{M}$ is conformally Osserman if and only if there exists a
quaternion structure $\{I,J,K\}$ on $V$ and constants $\lambda_I$, $\lambda_J$, $\lambda_K$ with $\lambda_I+\lambda_J+\lambda_K=0$
so that $R=\lambda_IR_I+\lambda_JR_J+\lambda_KR_K$ where $R_I$, $R_J$, and $R_K$ are given by Example {\rm \ref{exm-1} (2)}.
\end{enumerate}\end{theorem}

\section{Walker geometry}\label{sect-6}
One says $\mathcal{M}$ is a {\it Walker manifold} of signature $(2,2)$ if it admits a parallel totally
isotropic $2$-plane f\/ield; this implies \cite{Wa50,Wa50a}
that locally $\mathcal{M}$ is isometric to a metric on $\mathbb{R}^4$ with non-zero components
\begin{alignat*}{3}
& g(\partial_{x_1},\partial_{x_3})=g(\partial_{x_2},\partial_{x_4})=1,\qquad &&g(\partial_{x_3},\partial_{x_4})=g_{34},&\\
&g(\partial_{x_3},\partial_{x_3})=g_{33},&&g(\partial_{x_4},\partial_{x_4})=g_{44}.&
\end{alignat*}

The geometry of Walker manifolds with $g_{34}=0$ has been studied in \cite{CaGaMa05}. We impose a~dif\/ferent condition by setting
$g_{33}=g_{44}=0$ so the non-zero components of the metric are given~by:
\begin{gather}\label{eqn-6.a}
g(\partial_{x_1},\partial_{x_3})=g(\partial_{x_2},\partial_{x_4})=1\qquad\text{and}\qquad g(\partial_{x_3},\partial_{x_4})=g_{34}.
\end{gather}

By Theorem \ref{thm-11}, $\mathcal{M}$ is conformally Osserman if and only if $\cal{M}$ is either self-dual or anti-self-dual. One has
\cite{BrGaVa06} that:
\begin{theorem}\label{thm-12}
Let $\mathcal{M}=(\mathbb{R}^4,g)$ where $g$ is given by equation \eqref{eqn-6.a}.
\begin{enumerate}\itemsep=0pt
\item[\rm 1.] $\mathcal{M}$ is self-dual if and only if
$g_{34}=x_1p(x_3,x_4)+x_2q(x_3,x_4)+s(x_3,x_4)$.
\item[\rm 2.] $\mathcal{M}$ is anti-self-dual if and only if
$g_{34}=x_1p(x_3,x_4)+x_2q(x_3,x_4)+s(x_3,x_4) +\xi(x_1,x_4) +\eta(x_2,x_3)$ with
$p_{/3}=q_{/4}$ and  $g_{34}p_{/3}-x_1p_{/34}-x_2p_{/33}-s_{/34}=0$.
\end{enumerate}
\end{theorem}

We refer to \cite{BrGaGiVa07} for the following results:

\begin{theorem}\label{thm-13}
Let $\mathcal{M}=(\mathbb{R}^4,g)$ where $g$ is given by equation \eqref{eqn-6.a}. \begin{enumerate}\itemsep=0pt
\item[\rm 1.] The following conditions are equivalent:
\begin{enumerate}\itemsep=0pt
\item[\rm (a)] $\mathcal{M}$ is Osserman.
\item[\rm (b)] $\mathcal{M}$ is Einstein.
\item[\rm (c)] $\rho=0$.
\item[\rm (d)] $g_{34}=x_1p(x_3,x_4)+x_2q(x_3,x_4)+s(x_3,x_4)$ where
$p=-2a_4(a_0+a_3x_3+a_4x_4)^{-1}$, and
$q=-2a_3(a_0+a_3x_3+a_4x_4)^{-1}$ for $(a_0,a_3,a_4)\ne(0,0,0)$.
\item[\rm (e)] $\mathcal{J}(x)^2=0$ for all $x$.
\item[\rm (f)] $\mathcal{M}$ is Jacobi--Tsankov.\end{enumerate}
\item[\rm 2.] The following conditions are equivalent:\begin{enumerate}\itemsep=0pt
\item[\rm (a)] $\mathcal{M}$ is Jacobi--Videv.
\item[\rm (b)] $\mathcal{M}$ is skew-Tsankov.
\item[\rm (c)] $g_{34}=x_1p(x_3,x_4)+x_2q(x_3,x_4)+s(x_3,x_4)$ where
$p_{/3}=q_{/4}$.
\end{enumerate}\end{enumerate}\end{theorem}

A feature of these examples is that the warping functions are af\/f\/ine functions of $x_1$ and $x_2$. We return to
the general setting of Walker signature $(2,2)$ geometry.  Let $\nabla$ be a torsion free connection on a
$2$-dimensional manifold
$N$. Let
$(x_3,x_4)$ be local coordinates on $N$. We expand
\[
\nabla_{\partial_{x_i}}\partial_{x_j}=\sum_k\Gamma_{ij}{}^k\partial_{x_k}\qquad\text{for}\quad i,j,k=3,4
\]
 to
def\/ine the Christof\/fel symbols of $\nabla$. Let $\omega=x_1dx_3+x_2dx_4\in T^*N$; the pair $(x_1,x_2)$ gives
the dual f\/iber coordinates. Let $\xi=\xi_{ij}(x_3,x_4)\in C^\infty(S^2(T^*N))$ be an auxiliary
symmetric bilinear form.
\begin{definition}\label{defn-11}
The \emph{deformed Riemannian extension} is the Walker metric on $T^*N$
def\/ined by setting \cite{GaKuVa99}
\begin{gather*}
g(\partial_{x_1},\partial_{x_3})=g(\partial_{x_2},\partial_{x_4})=1,\\
g(\partial_{x_3},\partial_{x_3})=-2x_1\Gamma_{33}{}^3(x_3,x_4)
-2x_2\Gamma_{33}{}^4(x_3,x_4)+\xi_{33}(x_3,x_4),\\
g(\partial_{x_3},\partial_{x_4})=-2x_1\Gamma_{34}{}^3(x_3,x_4)
-2x_2\Gamma_{34}{}^4(x_3,x_4)+\xi_{34}(x_3,x_4),\\
g(\partial_{x_4},\partial_{x_4})=-2x_1\Gamma_{44}{}^3(x_3,x_4)
-2x_2\Gamma_{44}{}^4(x_3,x_4)+\xi_{44}(x_4,x_4).
\end{gather*}
\end{definition}

\begin{definition}\label{defn-12}
Let $\rho_N(x,y):=\operatorname{Tr}(z\rightarrow \mathcal{R}_\nabla(z,x)y)$ be the af\/f\/ine Ricci tensor. We may
decompose this
$2$-tensor into symmetric and anti-symmetric parts by def\/ining:
\begin{gather*}
\rho_N^s(x,y):=\tfrac12(\rho_N(x,y)+\rho_N(y,x))\qquad\text{and},\\
\rho_N^a(x,y):=\tfrac12(\rho_N(x,y)-\rho_N(y,x)).
\end{gather*}
The Jacobi operator is def\/ined by setting $\mathcal{J}_\nabla(x):y\rightarrow\mathcal{R}_\nabla(y,x)x$. We say that
$\mathcal{N}:=(N,\nabla)$ is {\it affine Osserman} if
$\mathcal{J}_\nabla(x)$ is nilpotent or, equivalently, if $\operatorname{Spec}\{\mathcal{J}_\nabla(x)\}=\{0\}$ for all $x$.
\end{definition}

We refer to \cite{BrGaGiVa07} for the proof of the following result:
\begin{theorem}\label{thm-14}
\ \begin{enumerate}\itemsep=0pt
\item[\rm 1.] $\mathcal{M}$ is skew-Tsankov if and only if $\rho_{N}^a=0$.
\item[\rm 2.] $\mathcal{M}$ is Osserman if and only if $\mathcal{N}$ is affine Osserman if and only if $\rho_N^s=0$.
\item[\rm 3.] $\rho_N^a=0$ or $\rho_N^s=0$ if and only if $\mathcal{M}$ is Jacobi--Videv.
\item[\rm 4.] $\rho_N=0$ if and only if $\mathcal{M}$ is Jacobi--Tsankov.
\end{enumerate}\end{theorem}

\begin{remark}\label{rmk-4}
This shows the notions Jacobi--Videv, and Jacobi--Tsankov, and skew-Tsankov are inequivalent notions.
\end{remark}

If $\mathcal{M}$ is conformally Osserman, let $m_\lambda$ be the
minimal polynomial of $\mathcal{J}_W$ and let $\operatorname{Spec}_W$ be the spectrum of $\mathcal{J}_W$. One has \cite{BrGaVa06}:
\begin{theorem}\label{thm-15}
 Let $\mathcal{M}=(\mathbb{R}^4,g)$ be the Walker manifold with non-zero metric components:
\[
g(\partial_{x_1},\partial_{x_3})=g(\partial_{x_2},\partial_{x_4})=1,\qquad\text{and}\qquad
g(\partial_{x_3},\partial_{x_4})=g_{34}.
\]
The following choices of
$g_{34}$ make $\mathcal{M}$ conformal Osserman with:
\begin{enumerate}\itemsep=0pt
\item[\rm 1.] The Jordan normal form does not change from point to point:
\begin{enumerate}\itemsep=0pt
\item[\rm (a)]  If $g_{34}=x_1^2-x_2^2$, then $m_\lambda=\lambda(\lambda^2-\frac14)$ and
$\operatorname{Spec}_W=\{0,0,\pm\frac12\}$.
\item[\rm (b)] If $g_{34}=x_1^2+x_2^2$, then $m_\lambda=\lambda(\lambda^2+\frac14)$ and
$\operatorname{Spec}_W=\big\{0,0,\pm\frac{\sqrt{-1}}2\big\}$.
\item[\rm (c)] If $g_{34}=x_1x_4+x_3x_4$, then $m_\lambda=\lambda^2$
and $\operatorname{Spec}_W=\{0\}$.
\item[\rm (d)]
If $g_{34}=x_1^2$, then $m_\lambda=\lambda^3$  and
$\operatorname{Spec}_W=\{0\}$.\end{enumerate}
\item[\rm 2.] $\operatorname{Spec}_W=\{0\}$ but the Jordan normal form changes from point to
point.
\begin{enumerate}\itemsep=0pt
\item[\rm (a)] If $g_{34}=x_2x_4^2+x_3^2x_4$, then $m_\lambda=\lambda^3$ if $x_4\ne0$,
  $m_\lambda=\lambda^2$ if $x_4=0$ and $x_3\ne0$, and $m_\lambda=\lambda$ if $x_3=x_4=0$.
\item[\rm (b)] If $g_{34}=x_2x_4^2+x_3x_4$, then $m_\lambda=\lambda^3$ if $x_4\ne0$, and $m_\lambda=\lambda^2$ if $x_4=0$.
\item[\rm (c)] If $g_{34}=x_1x_3^2$, then $m_\lambda =\lambda^3$ if $x_3\ne0$, and $m_\lambda=\lambda$ if $x_3=0$.
\item[\rm (d)] If $g_{34}=x_1x_3+x_2x_4$, then $m_\lambda =\lambda^2$ if $x_1x_3+x_2x_4\ne0$, and $m_\lambda=\lambda$ if $x_1x_3+x_2x_4=0$.
\end{enumerate}
\item[\rm 3.] The eigenvalues can change from point to point:
\begin{enumerate}\itemsep=0pt
\item[\rm (a)] If $g_{34}=x_1^4+x_1^2-x_2^4-x_2^2$, then
$\operatorname{Spec}_W=\big\{0,0,\pm{\textstyle\frac12}\sqrt{(6x_1^2+1)(6x_2^2+1)}\big\}$.
\item[\rm (b)] If $g_{34}=x_1^4+x_1^2+x_2^4+x_2^2$, then
$\operatorname{Spec}_W=\big\{0,0,\pm{\textstyle\frac12}\sqrt{-(6x_1^2+1)(6x_2^2+1)}\big\}$.
\item[\rm (c)] If $g_{34}=x_1^3-x_2^3$, then $\operatorname{Spec}_W=\{0,0,\pm{\textstyle\frac32}\sqrt{x_1x_2}\}$.
\end{enumerate}\end{enumerate}\end{theorem}

We conclude our discussion with the following result \cite{BrGiVa08}:

\begin{theorem}\label{thm-16}
Of the manifolds given above in Theorem {\rm \ref{thm-15}}, only the manifold with
$g_{34}=x_1^2$ is curvature homogeneous and only the manifold with $g_{34}=x_1x_4+x_3x_4$ is geodesically complete.
\end{theorem}

\subsection*{Acknowledgements} The research of M. Brozos-V\'azquez and of P. Gilkey
was partially supported by the Max Planck Institute for Mathematics in the Sciences (Germany)  and
by Project MTM2006-01432 (Spain). The research of E. Garc\'{\i}a--R\'{\i}o and of R. V\'azquez-Lorenzo was
partially supported by PGIDIT06PXIB207054PR (Spain). The research of S. Nik\v cevi\'c was partially supported by
Project 144032 (Serbia). It is a pleasure to acknowledge helpful conversations with C.~Dunn, E.~Puf\/f\/ini, and Z.~Zhelev concerning these and related matters.

\pdfbookmark[1]{References}{ref}

\LastPageEnding


\begin{thebibliography}{99}

\footnotesize\itemsep=0pt

\bibitem{BlGi04a} Bla{\v z}i{\'c}  N., Gilkey P.,
Conformally Osserman manifolds and conformally complex space
forms, {\it Int. J. Geom. Methods Mod. Phys.} \textbf{1}
(2004), 97--106, \href{http://arxiv.org/abs/math.DG/0311263}{math.DG/0311263}.

\bibitem{BlGi05}Bla{\v z}i{\'c}  N., Gilkey P., Conformally Osserman manifolds and self-duality in Riemannian
geometry, in  Proceedings of the Conference ``Dif\/ferential Geometry and Its
Applications'' (August 30~-- September 3, 2004, Charles University, Prague, Czech Republic), Editors J.~Bures,
O.~Kowalski, D.~Krupka and J.~Slovak, MATFYZPRESS, 2005, 15--18, \href{http://arxiv.org/abs/math.DG/0504498}{math.DG/0504498}.

\bibitem{BlGiNiSi05}
 Bla{\v z}i{\'c} N., Gilkey P., Nik\v cevi\'c S., Simon U.,
The spectral geometry of the Weyl conformal tensor,
\emph{Banach Center Publ.} \textbf{69} (2005), 195--203, \href{http://arxiv.org/abs/math.DG/0310226}{math.DG/0310226}.

\bibitem{BrGaGiVa07} Brozos-V\'azquez M., Garc\'{\i}a--R\'{\i}o E., Gilkey P., V\'azquez-Lorenzo R.,
Examples of signature $(2,2)$ manifolds with commuting curvature operators, \emph{J. Phys. A: Math. Theor.}, to appear, \href{http://arxiv.org/abs/0708.2770}{arXiv:0708.2770}.

\bibitem{BrGaVa06} Brozos-V\'azquez M., Garc\'{\i}a--R\'{\i}o E., V\'azquez-Lorenzo R., Conformally Osserman four-dimensional
manifolds whose conformal Jacobi operators have complex eigenvalues, \emph{Proc. Royal Soc. A} \textbf{462} (2006), 1425--1441.

\bibitem{BrGiVa08} Brozos-V\'azquez M., Garc\'{\i}a--R\'{\i}o E., Gilkey  P., V\'azquez--Lorenzo R.,
Completeness, Ricci blowup, the Osserman and the conformal
Osserman condition for Walker signature (2,2) manifolds,
in Proceedings of XV International Workshop on Geometry and
Physics,  to appear, \href{http://arxiv.org/abs/math.DG/0611279}{math.DG/0611279}.

\bibitem{BrGi06} Brozos-V\'azquez M., Gilkey P., Pseudo-Riemannian
manifolds with commuting Jacobi operators, \emph{Rend. Circ.
Mat. Palermo} \textbf{55} (2006), 163--174, \href{http://arxiv.org/abs/math.DG/0608707}{math.DG/0608707}.

\bibitem{BrGi07b} Brozos-V\'azquez M.,  Gilkey P.,
The global geometry of Riemannian manifolds with commuting
curvature operators, \emph{ J. Fixed Point Theory Appl.}
\textbf{1} (2007), 87--96, \href{http://arxiv.org/abs/math.DG/0609500}{math.DG/0609500}.

\bibitem{BrGi07} Brozos-V\'azquez  M., Gilkey P., Manifolds with
commuting Jacobi operators, \emph{J. Geom.} \textbf{86} (2007), 21--30, \href{http://arxiv.org/abs/math.DG/0507554}{math.DG/0507554}.

\bibitem{BrGiNi07} Brozos-V\'azquez M., Gilkey P., Nik\v cevi\'c S.,
Jacobi--Tsankov manifolds which are not 2-step nilpotent,
in Proceedings of the Conference ``Contemporary Geometry and Related Topics''
(June 26 -- July 2, 2005, Belgrade, Serbia and Montenegro), Editors N.~Bokan,  M.~Djori\'c,  A.T.~Fomenko, Z.~Rakic,
 B.~Wegner and J.~Wess,
University of Belgrade, Serbia,
2006, 63--79, \href{http://arxiv.org/abs/math.DG/0609565}{math.DG/0609565}.

\bibitem{CaGaMa05}Chaichi M., Garc\'{\i}a--R\'{\i}o E., Matsushita Y.,
Curvature properties of four-dimensional Walker metrics,
\emph{Classical Quantum Gravity} \textbf{22} (2005), 559--577.

\bibitem{DiGa04}
 D\'{\i}az-Ramos J.C., Garc\'{\i}a-R\'{\i}o E.,
A note on the structure of algebraic curvature tensors,
\emph{Linear Algebra Appl.} \textbf{382} (2004), 271--277.

\bibitem{Fi03} Fiedler B.,
Determination of the structure of algebraic curvature tensors by means of Young symmetrizers,
\emph{Seminaire Lotharingien de Combinatoire} \textbf{B48d} (2003), 20 pages,
\href{http://arxiv.org/abs/math.CO/0212278}{math.CO/0212278}.

\bibitem{FiGi03} Fiedler  B., Gilkey P.,
Nilpotent Szab\'o, Osserman and Ivanov--Petrova pseudo Riemannian
manifolds, \emph{Contemp. Math.} \textbf{337} (2003), 53--64, \href{http://arxiv.org/abs/math.DG/0211080}{math.DG/0211080}.

\bibitem{GaKuVa99}  Garc\'{\i}a--R\'{\i}o E., Kupeli D.N., V\'azquez-{Abal} M.E.,  V\'azquez-Lorenzo R.,
Af\/f\/ine Osserman connections and their Riemann extensions,
\emph{Differential Geom. Appl.} \textbf{11} (1999), 145--153.

\bibitem{GaKuVa02} Garc\'{\i}a-R\'{\i}o E., Kupeli D., V\'azquez-Lorenzo R.,
Osserman manifolds in semi-Riemannian geometry,
{\it Lecture Notes in Mathematics}, Vol.~1777,
Springer-Verlag, Berlin, 2002.

\bibitem{Gi01} Gilkey P.,
Geometric properties of natural operators def\/ined by the
Riemann curvature tensor, World Scien\-tif\/ic, 2001.

\bibitem{Gi07a} Gilkey P.,
The geometry of curvature homogeneous pseudo Riemannian manifolds,
Imperial College Press, 2007.

\bibitem{GiPuVi07} Gilkey P., Puf\/f\/ini E., Videv V., Puf\/f\/ini--Videv models and manifolds,  \emph{J. Geom.}, to
appear, \href{http://arxiv.org/abs/math.DG/0605464}{math.DG/0605464}.

\bibitem{GiNi07} Gilkey P., Nik\v cevi\'c S.,
Pseudo-Riemannian Jacobi--Videv manifolds, {\it Int. J. Geom. Methods Mod. Phys.} {\bf 4} (2007), 727--738,
\href{http://arxiv.org/abs/0708.1096}{arXiv:0708.1096}.

\bibitem{GiStVi98} Gilkey P., Stanilov G., Videv V.,
Pseudo Riemannian manifolds whose generalized Jacobi operator has constant characteristic polynomial,
\emph{J. Geom.} {\bf 62} (1998) 144--153.

\bibitem{IvViZh07} Ivanova M., Videv V., Zhelev Z.,
Four-dimensional Riemannian manifolds with commuting higher order Jacobi operators,
\href{http://arxiv.org/abs/math.DG/0701090}{math.DG/0701090}.

\bibitem{Os80} Osserman R.,
Curvature in the eighties,
\emph{Amer. Math. Monthly} \textbf{97} (1990), 731--756.


\bibitem{StVi92}
Stanilov G., Videv V.,
On a generalization of the Jacobi operator in the Riemannian geometry,
\emph{God. Sofij. Univ., Fak. Mat. Inform.} \textbf{86} (1994) 27--34.

\bibitem{StVi98} Stanilov  G., Videv V.,
Four dimensional pointwise Osserman manifolds,
\emph{Abh. Math. Sem. Univ. Hamburg} \textbf{68} (1998), 1--6.

\bibitem{StVi04}
 Stanilov G., Videv V.,
On the commuting of curvature operators,
in Proceedings  of the 33rd
Spring Conference of the Union of Bulgarian Mathematicians Borovtes ``Mathematics and Education in Mathematics''
(April 1--4, 2004, Sof\/ia), Sof\/ia, 2004, 176--179.

\bibitem{Ts05} Tsankov Y.,
A characterization of $n$-dimensional hypersurface in Euclidean space with commuting curvature operators,
\emph{Banach Center Publ.} \textbf{69} (2005), 205--209.

\bibitem{Vixx} Videv V., A characterization of the $4$-dimensional Einstein Riemannian manifolds using curvature
operators, Preprint.

\bibitem{Wa50} Walker A.G.,
Canonical form for a Riemannian space with a parallel f\/ield of null planes,
\emph{Quart. J. Math., Oxford Ser. (2)} \textbf{1}  (1950), 69--79.

\bibitem{Wa50a} Walker A.G.,
Canonical forms. II. Parallel partially null planes,
\emph{Quart. J. Math., Oxford Ser. (2)} \textbf{1} (1950), 147--152.

\end{thebibliography}
\end{document}